\documentclass[12pt]{article}

\usepackage{amsfonts}
\usepackage{}

\usepackage{amssymb}
\usepackage{amsmath}
\usepackage{graphics}
\usepackage{epsfig}
\newtheorem{claim}{\bf \t}[part]

\newtheorem{Definition}{Definition}[part]

\newtheorem{Lemma}{Lemma}[part]
\newtheorem{Proposition}{Proposition}[part]
\newtheorem{Remark}{Remark}[part]
\newtheorem{Theorem}{Theorem}[part]

\numberwithin{Assumption}{section} \numberwithin{Corollary}{section}
\numberwithin{Definition}{section} \numberwithin{equation}{section}
\numberwithin{Example}{section} \numberwithin{Lemma}{section}
\numberwithin{Proposition}{section} \numberwithin{Remark}{section}
\numberwithin{Theorem}{section}

\def\t{\theta}

\def\g{\gamma}

\def\f{\frac}

\def\text#1{{\rm #1}}

\begin{document}
\date{}
\title{\Large \bf Remarks on Blow-up of Smooth Solutions to the Compressible Fluid with Constant and Degenerate Viscosities}
\author{\small \textbf{Quansen Jiu},$^{1}$\thanks{The research is partially
supported by National Natural Sciences Foundation of China (No.
11171229, No.11231006 and No.11228102) and Project of Beijing Chang
Cheng Xue Zhe. E-mail: jiuqs@mail.cnu.edu.cn}\quad
  \textbf{Yuexun Wang}$^{2}$\thanks{ E-mail:
yx-wang13@mails.tsinghua.edu.cn}\quad
and \textbf{Zhouping Xin}$^{3}$\thanks{The research is partially
supported by Zheng Ge Ru Funds, Hong Kong RGC Earmarked Research
Grant CUHK4042/08P and CUHK4041/11P, and a grant from the Croucher
Foundation. Email: zpxin@ims.cuhk.edu.hk}} \maketitle \small $^1$
School of Mathematical Sciences, Capital Normal University, Beijing
100048, P. R. China

\small $^2$
Mathematical Sciences center, Tsinghua University, Beijing
100084, P. R. China

\small $^3$The Institute of Mathematical Sciences, Chinese
University of HongKong, HongKong\\

\small \maketitle { \bf Abstract:}  In this paper, we will show the
blow-up of smooth solutions to the Cauchy problem for the full
compressible Navier-Stokes equations and isentropic compressible
Navier-Stokes equations  with constant and degenerate viscosities in
arbitrary dimensions under some restrictions on the initial data. In
particular, the results hold true for the full compressible Euler
equations and isentropic compressible Euler equations and the
blow-up time can be computed in a more precise way.  It is  not
required that the initial data has compact support or contain vacuum
in any finite regions. Moreover, a simplified and unified proof on
the blow-up results to the classical solutions of the full
compressible Navier-Stokes equations without heat conduction by Xin
\cite{Xin} and with heat conduction by Cho-Bin \cite{CJ} will be
given.

\section{Introduction } \setcounter{equation}{0}
\setcounter{Assumption}{0} \setcounter{Theorem}{0}
\setcounter{Proposition}{0} \setcounter{Corollary}{0}
\setcounter{Lemma}{0}

    \qquad  The full compressible Navier-Stokes equations with constant viscosities
    read as
\begin{eqnarray}\label{CNS-1}
\left\{ \begin{array}{ll}
\partial_t\rho+\textrm{div}(\rho
u)=0,\\
\partial_t(\rho u)+\textrm{div}(\rho u \otimes u)+\nabla p=\textrm{div}(\mathcal{T}),\\
\partial_t({\f12}\rho |u|^2+\rho e )+\textrm{div}(({\f12}\rho |u|^2+\rho e+p)u)=\textrm{div}(u\mathcal {T})+\mathcal {K}
\Delta\theta.
\end{array}
\right.
\end{eqnarray}
Here $(x,t)\in\mathbb{R}^n\times\mathbb{R}_+ $ and $\rho= \rho
(x,t), u = (u_1, u_2,\cdots ,u_n),\theta, p $ and $e$ denote the
density, velocity, absolute temperature, pressure and internal
energy, respectively. $\mathcal {T}$ is the stress tensor given by
$$\mathcal {T}=\mu(\nabla u+(\nabla u)^t)+\lambda \textrm{div}uI,$$
where $I$ is the identity matrix, and $\mu$ and $\lambda$ are the
coefficients of viscosity and the second coefficient of viscosity,
respectively, which satisfy
$$\mu\geq0,2\mu+n\lambda\geq0.$$
 We also denote by $ \mathcal {K}\geq0 $
the coefficient of heat conduction.

If $\mu= \lambda=  \mathcal {K}= 0$, the Navier-Stokes equations
\eqref{CNS-1} become the
 compressible Euler equations.  The polytropic gas satisfies the following state  equations:
\begin{equation}
p = R\rho\theta,\ e = c_\nu\theta, \quad p =
Aexp({\f{s}{c_\nu}})\rho^\gamma\label{(4)},
\end{equation}
where  $R > 0$ is the gas constant, $A>0$ is an absolute constant,
$\gamma> 1$  is the specific heat ratio, $c_\nu= {\f{R}{\gamma-1}}$
and $s$ is the entropy. For simplicity, we take $A=1$. The pressure
can be expressed as
\begin{equation}
p=(\gamma-1)\rho e.
\end{equation}

The initial data to  the equations \eqref{CNS-1} are imposed as
 \begin{equation}\label{6}
  (\rho,u,s)(x,t)|_{t=0}=(\rho_0(x),u_0(x),s_0(x)).
 \end{equation}

If the entropy is a constant, the full compressible Navier-Stokes
equations reduce to isentropic ones which read as
\begin{eqnarray}\label{CNS-2}
\left\{ \begin{array}{ll}
\partial_t\rho+\textrm{div}(\rho
u)=0,\\
\partial_t(\rho u)+\textrm{div}(\rho u \otimes u)+\nabla p=\textrm{div}(\mathcal
{T}).\\
\end{array}
\right.
\end{eqnarray}

The state  equation of the isentropic process becomes
\begin{equation}\label{5}
  p=\rho^\gamma, \gamma>1.
 \end{equation}

The initial data to the equations \eqref{CNS-2} are imposed as

\begin{equation}\label{7}
  (\rho,u)(x,t)|_{t=0}=(\rho_0(x),u_0(x)).
 \end{equation}

 The compressible Navier-Stokes equations with density-dependent viscosity are another kind of important models.
  When deriving the compressible Navier-Stokes equations by
Chapman-Enskog expansions from the Boltzmann equation, the viscosity
 depends on the temperature and thus on the density for the isentropic
flows. The isentropic compressible Navier-Stokes equations with
density-dependent viscosity can be written as
\begin{eqnarray}\label{CNS-3}
\left\{ \begin{array}{ll}
\partial_t\rho+\textrm{div}(\rho u)=0,\\
\partial_t(\rho u)+\textrm{div}(\rho u \otimes u)+\nabla p=\textrm{div}[h(\rho)({\f{\nabla u+(\nabla u)^t}{2}})]+\nabla
(g(\rho)\textrm{div}u),\\
\end{array}
\right.
\end{eqnarray}
where the pressure $p$ is same as in \eqref{5}, $h(\rho)$ and
$g(\rho)$ are the Lam\'{e} viscosity coefficients satisfying
\begin{eqnarray} h(\rho)\ge 0,
h(\rho)+ng(\rho)\geq 0.\label{LVC}
\end{eqnarray}

In particular, we consider the case
$$
h(\rho)=\rho^\alpha, g(\rho)=(\alpha-1)\rho^\alpha $$
for$\alpha>1-\frac1n$ such that \eqref{LVC} is satisfied. For more
general $h(\rho)$ and $g(\rho)=\rho h'(\rho)-h(\rho)$, it is
referred to \cite{BDL}, \cite{GJX} and references therein for
 entropy estimates and existence of weak solutions.  It should be
noted that the viscous Saint-Venant system for the shallow water,
derived from the incompressible Navier-Stokes equations with a
moving free surface, is expressed exactly as in \eqref{CNS-3} with
$h(\rho)=\rho, g(\rho)=0$ and $\g=2$.

There have been extensive studies on the compressible Navier-Stokes
equations with constant and degenerate viscosities (see
\cite{CJ},\cite{CK1},\cite{CK2},\cite{FE1},
\cite{HLX4},\cite{HD},\cite{HSD},\cite{HSJ},
 \cite{PLL},\cite{MN1},\cite{MN2}
 on the case of  constant viscosity
and
\cite{DJ},\cite{JWX},\cite{JWX-2},\cite{JX},\cite{LLX},\cite{MV2},
\cite{Kazhikhov},\cite{JWX-2} on the case of degenerate viscosity).
For the Navier-Stokes equations with constant viscosity, the global
existence and uniqueness of the strong (classical)
 solution in  one-dimensional case has been well-understood (see \cite{Ka},
 \cite{KS1977}, \cite{HSJ} and references therein)
and it is proved that if the solution has no vacuum initially then
it will not appear vacuum later (in any finite time). Moreover, even
if vacuum is permitted, the global well-posedness of the solution to
the one-dimensional navier-Stokes equations was studied recently
under some compatibility conditions (see \cite{DWZ},\cite{JLY-2013}
and references therein). However, in multi-dimensional case, the
global well-posedness of the classical solution to these models for
large initial data remains completely open. For the Navier-Stokes
equations with density-dependent viscosity, the global existence and
uniqueness of the classical solution remain open in
multi-dimensional case except the periodic or Cauchy problem of the
2D Kazhikhov-Vaigant model in compressible flow (see
\cite{Kazhikhov}, \cite{JWX-2} and references therein). In fact, if
vacuum is permitted,  the global well-posedness of the classical
solution to the one-dimensional Navier-Stokes equations with
density-dependent viscosity still remains open.

These will be  subtle issues because the classical solution to the
compressible Euler equations and Navier-Stokes equations may blow up
in general. For the compressible Euler equations, Sideris
\cite{Sideris} firstly showed that the life span of the classical
solution is finite if the initial velocity is large enough in some
region with compact support. Makino,  Ukai and Kawashima \cite{MUK}
studied the blow-up of the classical solution if initial density and
velocity hold compact support. Chemin \cite{Chemin} investigated the
blow-up of the classical solution to the one-dimensional
compressible Euler equations. Recently, the 3-D shock formation in
general settings was studied by Christodoulou \cite{Christodoulou}.
For the compressible Navier-Stokes equations without heat
conduction, Xin \cite{Xin} firstly obtained the blow-up of the
classical solution under assumptions that the initial density has
compact support, which was generalized by Cho and Jin to the case
with heat conduction in \cite{CJ} and Razanova
 to the case with initial data rapidly decay at far fields in
\cite{Rozanova}. Luo and Xin \cite{LX} studied the blow-up of
symmetric smooth solutions to two dimensional isentropic
Navier-Stokes equations. Recently,
 Xin and Yan  \cite{XY} proved that any
classical solutions of viscous compressible fluids without heat
conduction will blow up in finite time, as long as the initial data
has an isolated mass group.

In this paper, we will prove the blow-up of the classical solution
to the Cauchy problem for compressible Navier-Stokes equations  in
arbitrary dimensions with constant viscosity and degenerate
viscosities under some restrictions of the initial data.  In
particular, the results hold true for the compressible Euler
equations  and the blow-up time can be computed in a more precise
way.  It should be noted that it is not required that the initial
data has compact support or has fast decay at far fields. And we
will give a new  and simplified proof of the blow-up results
appeared in \cite{Xin} and \cite{CJ}. To obtain our main results,
some physical quantities such as mass, momentum, momentum of
inertia, internal energy, potential energy, total energy and some
combined functionals of these quantities are introduced. The basic
properties and specific  relationships between them are helpful and
crucial to prove the main results.  More precisely, for the full
Navier-Stokes and Euler equations, the upper and lower decay rates
of the internal energy will be calculated in a precise way. And for
the isentropic Navier-Stokes and Euler equations, the upper and
lower decay rates of the potential energy will be presented
accordingly. Then the main results will be proved by comparing the
coefficients of the upper and lower decay rates.

In the following, we denote the full Navier-Stokes equations
\eqref{CNS-1} and \eqref{(4)} as (CNS), the isentropic Navier-Stokes
equations \eqref{CNS-2} and \eqref{5} as (ICNS) and the Naver-Stokes
equations with density dependent viscosity  \eqref{CNS-3} and
\eqref{5} as (DICNS) for short.

The paper is organized as follows. In section 2, we will introduce
some physical quantities and present our main results. In Section 3
and 4, we will give some basic properties of the physical quantities
and the proof of the main results.

\section{Main Results }

The following  physical quantities will be  used in this paper.
\begin{equation}\label{Q1}
\left\{
\begin{array}{lll}
&M(t)=\int _{\mathbb{R}^n}\rho  dx \ \ &(mass),\\
&\mathbb{P}(t)=\int_{\mathbb{R}^n} \rho u dx\ \ &(momentum),\\
&F(t)=\int_{\mathbb{R}^n} \rho u x dx &(momentum\   weight), \\
&G(t)={\f12}\int_{\mathbb{R}^n} \rho |x|^2 dx \ \ &(momentum \ of \ inertia),\\
&E(t)={\f12}\int _{\mathbb{R}^n} \rho |u|^2 dx+\int_{\mathbb{R}^n}
\rho e dx\triangleq
E_k(t)+E_i(t)\ \ &(total \ energy),\\
&IE(t) ={\f12}\int_{\mathbb{R}^n} \rho |u|^2 dx+{\f{1}{\gamma
-1}}\int_{\mathbb{R}^n} P dx\triangleq E_k(t)+I(t) \ \ &(energy).
\end{array}
\right.
\end{equation}

The basic properties and relationships between these quantities will
be discussed in Section 3 and 4. To prove our main results, we
introduce the following functionals:
\begin{equation}\label{Q2}
\left\{
\begin{array}{ll}
&H(t)=2 E_k(t)+n(\gamma -1)E_i(t),\\
&IH(t)=2 E_k(t)+n(\gamma
-1)I(t),\\
&J(t)= G(t)-(t+1)F(t)+(t+1)^2 E(t),\\
& IJ(t)= G(t)-(t+1)F(t)+(t+1)^2 IE(t).
\end{array}
 \right.
\end{equation}

We always assume that $M(0), \mathbb{P}(0), F(0), G(0),  E(0),
IE(0)$ are finite and   $M(0)>0, \mathbb{P}(0)\neq 0,  E(0)>0,
IE(0)>0$. Since $F(t)^2\leq4G(t)E_k(t)$ (see Lemma \ref{(101)}) and
$E_i(t)\ge 0, I(t)\ge 0$, it follows that $J(0)>0$ and $IJ(0)>0$
respectively.

For technical reason, we impose decay conditions on the solutions as
follows:
\begin{eqnarray}\label{9-25}
|u|\rightarrow0,\rho u=o({{\f{1}{|x|^{n+1}}}}),|\nabla
u|=o({{\f{1}{|x|^n}}}),P=o({{\f{1}{|x|^n}}}),
\end{eqnarray}
as $|x|\to \infty$ for any fixed $t>0$.

If $\mathcal {K}\neq 0$, we impose
\begin{eqnarray}\label{9-25-1}
|\nabla \theta|=o(\f{1}{|x|^{n-1}}), \ |x|\to\infty, \ t>0.
\end{eqnarray}

It should be remarked that the conditions
\eqref{9-25}-\eqref{9-25-1}  guarantee that the integration by parts
in our calculations make sense, which is similar to those in
\cite{Rozanova}. In particular, for the classical solution to (CNS)
satisfying \eqref{9-25}-\eqref{9-25-1}, the conservations of the
mass $M(t)$, momentum $\mathbb{P}(t)$ and energy $E(t)$ hold true.

 Before we state our main results, we give the following definitions.

\begin{Definition}\label{(0000001)}
For the Cauchy problem  of (CNS) without heat conduction ($\mathcal
{K}=0$), we call $ (\rho,u,s)\in X(T)$ if $ (\rho,u,s)$ is a
classical solution in $[0,T]$ satisfying \eqref{9-25}.
\end{Definition}

\begin{Definition}\label{(0000002)}
For the Cauchy problem of (ICNS), we call $ (\rho,u)\in Y(T)$ if $
(\rho,u)$ is a classical solution in $[0,T]$ satisfying
\eqref{9-25}.
\end{Definition}

\begin{Definition}\label{(0000003)}
For the Cauchy problem of  (DICNS), we call $ (\rho,u)\in Z(T)$ if $
(\rho,u)$ is a classical solution in $[0,T]$ satisfying \eqref{9-25}
in which the condition $|\nabla u|=o({{\f{1}{|x|^n}}})$ is replaced
by
\begin{eqnarray}\label{9-25-2}
\rho^\alpha|\nabla u|=o({{\f{1}{|x|^n}}}).
\end{eqnarray}
\end{Definition}

Our main results  are stated as follows.
 \begin{Theorem}\label{(1000)}
  Let $1<\gamma \leq1+{\f{2}{n }}$. If  the initial values satisfy
\begin{equation}\label{1000-1}
{\f {E(0)^{\f {n(\gamma-1)}{2}}J(0)}{\exp ({\f{s_1}{c_\nu}})M(0)^{\f
{(n+2)\gamma-n}{2}}}}<({\f{\Gamma ({\f n 2}+1)}{(\pi)^{\f n
2}}})^{\gamma-1 } {\f{1}{2^{\f{(n+2)\gamma-n}{2}}(\gamma-1)}},
\end{equation}
       where $s_1=\min_x s_0(x)$. Then there  exists a $T_1^{*}>0$ such that  there is no solution in $X(T_1^{*})$ to the
Cauchy problem of (CNS)  without heat conducting. Moreover, the
result holds true for the compressible Euler equations.
\end{Theorem}
\begin{Theorem}\label{(1001)}
Let $1<\gamma \leq1+{\f{2}{n }}$. If the initial values satisfy
 \begin{equation}\label{1000-2}
      {\f {E(0)^{\f {n(\gamma-1)}{2}}IJ(0)}{M(0)^{\f {(n+2)\gamma-n}{2}}}}<({\f{\Gamma ({\f n 2}+1)}{(\pi)^{\f n 2}}})^{\gamma-1 }{\f {1}{2^{\f{(n+2)\gamma-n}{2}}(\gamma-1)}}.
\end{equation}
Then there   exists a $T_2^{*}>0$ such that  there is no solution in
$Y(T_2^{*})$ to the Cauchy problem of (ICNS). In particular, the
claim above is true for isentropic compressible Euler equations.
\end{Theorem}

To study the blow-up of the classical solution to the isentropic
Navier-Stokes equations with density-dependent viscosity, we assume
that the energy in general n-dimensional case and the upper bound of
the density in 1-dimensional case are finite.

\vspace{3mm}

{\bf Assumptions 1. Energy Bounds:}
  \begin{eqnarray}
 \nonumber && \int_{\mathbb{R}^n} [{\f12}\rho u^2+{\f{1}{\gamma-1}}\rho^\gamma]dx+\int_0^t\int_{\mathbb{R}^n} h(\rho)|{\f{\nabla u+(\nabla u)^t}{2}}|^2 dxdt+\int_0^t\int_{\mathbb{R}^n} g(\rho)\textrm({div}u)^2dxdt\\
  &=&\int_{\mathbb{R}^n} [{\f12}\rho_0 u_0^2+{\f{1}{\gamma-1}}\rho_0^\gamma]dx
  \equiv C_{13}<\infty.\label{(43.1)}
  \end{eqnarray}

\vspace{3mm}

{\bf Assumptions 2. Upper Bound of  Density:}  When $n=1$, we assume
that $\max \rho(x,t)= C_{14}<\infty$.

\vspace{3mm}

Then we have
\begin{Theorem}\label{(1002)} Suppose that Assumptions 2 holds.
 Let $n=1$, $\alpha\geq\gamma$ and $1<\gamma \leq3$. If the initial values satisfy
  \begin{equation}\label{1000-3}
    \frac{IJ(0)\{{\f{1}{2}}[\max(2,(\gamma-1))IE(0)+C_{18}]\}^{\f{\gamma-1}{2}}}{M(0)^{\f
    {3\gamma-1}{2}}}< \f{\exp({\f{1-\gamma}{4}}C_{14}^{\alpha-\gamma})}{4^{2\gamma-1}(\gamma-1)}
\end{equation}
where  $C_{18}$ is any positive number between $0$ and
${\f{\mathbb{P}(0)^2}{M(0)}}$ and $C_{14}=\max\rho=C(\rho_0,u_0)$,
then there exists a $T_3^{*}>0$ such that  there is no solution in
$Z(T_3^{*})$ to the Cauchy problem of (DICNS).
\end{Theorem}
\begin{Theorem}\label{(1003)} Suppose that Assumptions 1 holds.
Let $n=1, {\f{\gamma+1}{2}}<\alpha\leq\gamma$ and $1<\gamma < 3$. If
the initial values satisfy
  \begin{eqnarray}\label{1000-4}
\f{IJ(0)
\{\f{1}{2}[\max(2,(\gamma-1))IE(0)+C_{18}]\}^{\f{\gamma-1}{2}}}{M(0)^\f{3\gamma-1}{2}\exp({\f{\alpha(1-\gamma)}{4(2\alpha-\gamma-1)}}(\gamma-1)^{\f{\alpha-1}{\gamma-1}}
M(0)^{\f{\gamma-\alpha}{\gamma-1}})}<\f{(\Gamma
(\f32))^{\gamma-1}}{(\gamma-1)\pi^{\gamma-1}2^{\frac{3\gamma-n}{2}}
},
\end{eqnarray}
where  $C_{18}$ is any positive constant between $0$ and
${\f{\mathbb{P}(0)^2}{M(0)}}$, then there exists a $T_4^{*}>0$ such
that there is no solutions
  in $Z(T_4^{*})$ to the Cauchy problem of (DICNS).
\end{Theorem}
\begin{Theorem}\label{(1003+)} Suppose that Assumptions 1 holds.
Let $n\ge 2, {\f{\gamma+1}{2}}<\alpha\leq\gamma$ and $1<\gamma
<1+{\f{2}{n }}$. If the initial values satisfy
  \begin{eqnarray}\label{1000-4+}
\frac{ IJ(0)
\{{\f{1}{2}}[\max(2,n(\gamma-1))IE(0)+C_{22}]\}^{\f{\gamma-1}{2}}}{
M(0)^{\f {(n+2)\gamma-n}{2}}\exp({\f{1-\gamma}{4(\alpha-1)
(2\alpha-\gamma-1)}}[1+n(\alpha-1)]^2
(\gamma-1)^{\f{\alpha-1}{\gamma-1}}M(0)^{\f{\gamma-\alpha}{\gamma-1}})}<\nonumber\\
\f{(\Gamma ({\f n 2}+1))^{\gamma-1}}{(\gamma-1)\pi^{\f n
2(\gamma-1)}2^{\f{(n+2)\gamma-n}{2}}},
\end{eqnarray}
where $C_{22}$ is any positive constant between $0$ and
${\f{\mathbb{P}(0)^2}{M(0)}}$, then there exists a $T_5^{*}>0$ such
that there is no solutions
  in $Z(T_5^{*})$ to the Cauchy problem  of (DICNS).
\end{Theorem}

A few remarks are in order.
\begin{Remark}
In comparison with results obtained by Rozanova  \cite{Rozanova},
the conditions imposed on the initial data are different from those
in \cite{Rozanova} and the blow-up of the classical solution to the
Navier-Stokes equations with density-dependent viscosities is
addressed here.
\end{Remark}
\begin{Remark}
For the one-dimensional Cauchy problem of (DICNS) with
$\mu(\rho)=\rho^\alpha$, Mellet-Vasseur \cite{MV2} proved the global
existence of  the strong solution $\rho\in L^\infty(0,T;
H^1(\mathbb{R}))$, $u\in L^\infty(0,T;$\\ $H^1(\mathbb{R}))\cap
L^2(0,T; H^2(\mathbb{R}))$ under the assumption $0<\alpha<\frac12$
and uniqueness under assumptions $ \mu(\rho)\ge \nu>0$ and
$\gamma\ge 2$ additionally. Under the jump free boundary conditions,
Jiang-Xin-Zhang \cite{JXZ-2005} established the global
well-posedness of the strong solutions for $0<\alpha<1$.  If
$\mu(\rho)=1+\rho^\beta$ with $\beta\ge 0$, Jiu-Li-Ye
\cite{JLY-2013} obtained the global well-posedness of the classical
solution permitting vacuum. However, Theorem \ref{(1002)} shows
that, if $\mu(\rho)=\rho^\alpha$ with $\alpha\ge \gamma, 1<\gamma\le
3$, the classical solution to the 1D Cauchy problem of (DICNS) will
blow up in general even though the density has upper bound. While if
$\mu(\rho)=\rho^\alpha$ with $\frac{\gamma+1}{2}\le \alpha< \gamma,
1<\gamma< 3$, Theorem \ref{(1003)} shows that the classical solution
to the 1D Cauchy problem of (DICNS) will blow up in general even
though the energy is bounded. The blow-up result in the
multi-dimensional case is presented in Theorem \ref{(1003+)}.
\end{Remark}
\begin{Remark}
In \eqref{1000-1}-\eqref{1000-4+}, we do not  require that the
initial data has compact support or contain  vacuum in any finite
region.
\end{Remark}
\begin{Remark}
It is clear that the equality \eqref{1000-1} holds as long as
initial entropy is sufficient large.
\end{Remark}
\begin{Remark}
 The time $T_1^{*}$-$T_5^{*}$ can be computed precisely, for example,  we can solve out
 $T_1^{*}$ in  Theorem \ref{(1001)}  by \eqref{T1}. In other words, we can
 find out the "last" blow-up time.
\end{Remark}
\begin{Remark}
  It should be noted that  Assumption 1 is the usual energy estimate and can be verified under the condition $\int_{\mathbb{R}^n} [{\f12}\rho_0
  u_0^2+{\f{1}{\gamma-1}}\rho_0^\gamma]dx<\infty$. The Assumptions 2 can also be verified in one-dimensional case under
  suitable conditions of the initial data( see \cite{GJ},\cite{JXZ-2005},\cite{JX},\cite{DJ} and references therein).
\end{Remark}
\begin{Remark}
  Our results show that if the conditions  in Theorem \ref{(1001)} are satisfied, the weak solutions obtained in
  \cite{FE1} and \cite{PLL} can not belong to  $Y(T_2^{*})$.
\end{Remark}

\section{ The Proof of Theorem \ref{(1000)}- \ref{(1001)}}

Let $(\rho,u,s)\in X(T) $ be a classical solution to the Cauchy
problem of (CNS). And $(\rho,u)\in Y(T) $ is a classical solution to
the Cauchy problem  of (ICNS). In this subsection, we will first
present some basic relationships between the quantities defined in
Section 2. Then we will give the proof of Theorem \ref{(1000)} and
Theorem \ref{(1001)}.
\begin{Lemma}\label{(100)} For (CNS) and (ICNS), we have
\begin{eqnarray}
{\f {d}{dt}}M(t)=0,  {\f {d}{dt}}\mathbb{P}(t)=0,  {\f {d}{dt}}G(t)=F(t).\label{(12)}
\end{eqnarray}
For (CNS), we have
\begin{eqnarray}
{\f {d}{dt}}E(t)=0,  {\f {d}{dt}}F(t)=H(t).\label{(13)}
\end{eqnarray}
For (ICNS), we have
\begin{eqnarray}
{\f {d}{dt}}IE(t)&&=-\int_{\mathbb{R}^n}
[2\mu\sum_{j=1}^n(\partial_ju_j)^2+\lambda(
\textrm{div}u)^2\nonumber\\
&&+\mu\sum_{i\neq
j}(\partial_ju_i)^2+2\mu\sum_{i>j}(\partial_ju_i)(\partial_iu_j)]dx,\label{(14+)}
\end{eqnarray}
\begin{eqnarray}
{\f {d}{dt}}F(t)=IH(t)&&.\label{(14)}
\end{eqnarray}
\end{Lemma}
 \textbf{{\em Proof.}} Using (CNS) and (ICNS), applying integration by
 parts,  one can verify \eqref{(12)}-\eqref{(14)}.

\begin{Lemma}\label{(101)}
 For(CNS) and(ICNS), we have
\begin{equation}
F(t)^2\leq4G(t)E_k(t),\label{(15)}
\end{equation}
and
\begin{equation}
\mathbb{P}(0)^2\leq2M(0)E_k(t).\label{(16)}
\end{equation}
\end{Lemma}
\textbf{{\em Proof.}} Using Lemma \ref{(100)} and $H\ddot{o}lder's$
inequality, one can verify \eqref{(15)}-\eqref{(16)}.

Based on Lemma \ref{(100)} and Lemma \ref{(101)}, we can obtain
two-sided estimates of G(t) (see \cite{Rozanova}).

\begin{Lemma}\label{(102)}
 For (CNS), if $1<\gamma \leq1+{\f{2}{n }}$, we have

\begin{equation}
{\f{n(\gamma-1)}{2}}E(0)t^2+F(0)t+G(0)\leq G(t)\leq E(0)t^2+F(0)t+G(0)\label{(17)},
\end{equation}
if $\gamma >1+{\f{2}{n }}$, we have
\begin{equation}
E(0)t^2+F(0)t+G(0)\leq G(t)\leq {\f{n(\gamma-1)}{2}}E(0)t^2+F(0)t+G(0)\label{(18)}.
\end{equation}
For (ICNS), if $1<\gamma \leq1+{\f{2}{n }}$, we have

\begin{equation}
{\f{\mathbb{P}(0)^2}{2M(0)}}t^2+F(0)t+G(0)\leq G(t)\leq
IE(0)t^2+F(0)t+G(0)\label{(19)},
\end{equation}
if $\gamma >1+{\f{2}{n }}$, we have

\begin{equation}
{\f{\mathbb{P}(0)^2}{2M(0)}}t^2+F(0)t+G(0)\leq G(t)\leq
{\f{n(\gamma-1)}{2}}IE(0)t^2+F(0)t+G(0)\label{(20)}.
\end{equation}
\end{Lemma}
\textbf{{\em Proof.}} For (CNS), in view of Lemma \ref{(100)}, if
$1<\gamma \leq1+{\f{2}{n }}$, we have

\begin{eqnarray}
     {\f {d^2}{dt^2}}G(t)\nonumber
     &=&H(t)
     =2E(t)+(n(\gamma-1)-2)E_i(t)\\
     &\leq&2E(t)
     =2E(0),\label{(21)}
\end{eqnarray}
and

\begin{eqnarray}
   {\f {d^2}{dt^2}}G(t)\nonumber
   &=&H(t)
   =n(\gamma-1)E(t)+(2-n(\gamma-1))E_k(t)\\
   &\geq& n(\gamma-1)E(t)
   =n(\gamma-1)E(0).\label{(22)}
\end{eqnarray}
 Integrating (\ref{(21)}) and (\ref{(22)}) over $[0,t]$, we get  (\ref{(17)}). The proof of \eqref{(18)} is similar.

For (ICNS), in view of Lemma \ref{(100)} and Lemma \ref{(101)}, if
$1<\gamma \leq1+{\f{2}{n }}$, we have

\begin{eqnarray}
     \nonumber{\f {d^2}{dt^2}}G(t)
     &=&IH(t)
     =2IE(t)+(n(\gamma-1)-2)I(t)\\
     &\leq&2IE(t)
     \leq2IE(0),\label{(23)}
\end{eqnarray}
and
\begin{eqnarray}
     {\f {d^2}{dt^2}}G(t)=IH(t)
   \geq2E_k(t)
     \geq{\f{\mathbb{P}(0)^2}{M(0)}},\label{(24)}
\end{eqnarray}
where \eqref{(16)} has been used. Integrating (\ref{(23)}) and
(\ref{(24)}) over $[0,t]$, we get (\ref{(19)}). The proof of
\eqref{(20)} is similar. We end the proof of the lemma.

The following lemma is due to Chemin  \cite{Chemin}.

\begin{Lemma}\label{(103)} For any $f\in L^1(\mathbb{R}^n,dx)\cap L^\gamma(\mathbb{R}^n,dx)\cap
L^1(\mathbb{R}^n,|x|^2dx)$, it holds that
\begin{eqnarray}
   \parallel f \parallel_{L^1(\mathbb{R}^n,dx)}\leq C_1\parallel f \parallel_{L^\gamma(\mathbb{R}^n,dx)}^{\f{2\gamma}{(n+2)\gamma-n}}\parallel f \parallel_{L^1(\mathbb{R}^n,|x| ^2dx)}^{\f{n(\gamma-1)}{(n+2)\gamma-n}},\label{(25-)}
\end{eqnarray}
where $C_1=2|B_1|^{\f{2(\gamma-1)}{(n+2)\gamma-n}}=2({\f{\pi^{\f n 2}}{\Gamma ({\f n 2}+1)}})^{\f{2(\gamma-1)}{(n+2)\gamma-n}}$.
\end{Lemma}
 \textbf{{\em
Proof.}} For any $r>0$, it follows from the $H\ddot{o}lder's$ inequality that

\begin{eqnarray}
  \int_{\mathbb{R}^n} |f(x)|dx&=&\int_ {|x|\leq r} |f(x)|dx+\int_ {|x|\geq r} |f(x)|dx\\
  &\leq& |B_r|^{1-{\f 1 \gamma}}(\int_ {|x|\leq r} |f(x)|^\gamma dx)^{\f 1 \gamma}+r^{-2}\int_ {|x|\geq r} |f(x)||x|^2dx)\\
  &\leq&|B_1|^{1-{\f 1 \gamma}}r^{n(1-{\f 1 \gamma})}(\int_{\mathbb{R}^n} |f(x)|^\gamma dx)^{\f 1 \gamma}+r^{-2}
  \int_{\mathbb{R}^n}|f(x)||x|^2 dx.\label{(25)}
\end{eqnarray}
Choosing

$$r=(\f{\parallel f \parallel_{L^1(\mathbb{R}^n,|x| ^2dx)}}{\parallel f \parallel_{L^\gamma(\mathbb{R}^n,dx)}|B_1|^{1-{\f{1}{\gamma}}}})^{{\f{\gamma}{(n+2)\gamma-n}}},$$
we obtain \eqref{(25-)}. The proof of the lemma is finished.

Taking $f=\rho$ in Lemma \ref{(103)}, we arrive at the lower bound
of $E_i(t) $ and $I(t)$, which is

\begin{Proposition}\label{(104)}
 For (CNS) and(ICNS), we have

\begin{eqnarray}
E_i(t)\geq {\f {C_2}{G(t)^{\f {n (\gamma-1)}{2}}}},\label{(26)}
\end{eqnarray}
and

\begin{equation}
  I(t)\geq {\f {C_3}{G(t)^{\f {n (\gamma-1)}{2}}}},\label{(27)}
\end{equation}
respectively, where $C_2=({\f{\Gamma ({\f n 2}+1)}{(\pi)^{\f n
2}}})^{\gamma-1 } {\f{\exp ({\f {s_1}{c_\nu}})M(0)^{\f
{(n+2)\gamma-n}{2}}}{2^{\f{(n+2)\gamma-n}{2}}(\gamma-1)}}$ and
$C_3=({\f{\Gamma ({\f n 2}+1)}{(\pi)^{\f n 2}}})^{\gamma-1 }
{\f{M(0)^{\f
{(n+2)\gamma-n}{2}}}{2^{\f{(n+2)\gamma-n}{2}}(\gamma-1)}}$.
\end{Proposition}

The following lemma was shown by Xin  \cite{Xin}. Here we give a new
proof of it.

\begin{Lemma}\label{(105)}
 For (CNS) and the (ICNS), the following estimates hold:
\begin{eqnarray}
{\f {d}{dt}}J(t)\leq
\begin{cases}
   {\f{2-n(\gamma-1)}{t+1}}J(t),\qquad 1<\gamma \leq1+{\f{2}{n }},\cr
    0,\qquad\qquad\qquad \qquad \gamma >1+{\f{2}{n }}.
   \end{cases}\label{(28)}
\end{eqnarray}
and
\begin{eqnarray}
{\f {d}{dt}}I J(t)\leq
\begin{cases}
   {\f{2-n(\gamma-1)}{t+1}}I J(t),\qquad 1<\gamma \leq1+{\f{2}{n }},\cr
    0,\qquad\qquad\qquad \qquad \gamma >1+{\f{2}{n }},
   \end{cases}\label{(29)}
\end{eqnarray}
respectively.
\end{Lemma}
\textbf{{\em
Proof.}} Due to Lemma \ref{(101)}, if we regard
$$ G(t)-(t+1)F(t)+(t+1)^2 E_k(t)$$
as a quadratic function of $(t+1)$, since
$$\Delta=(F(t)^2-4G(t)E_k(t))\leq0,$$
 we have
$$G(t)-(t+1)F(t)+(t+1)^2 E_k(t)\geq0.$$
Consequently,
\begin{eqnarray}
  E_i(t)\leq {\f{1}{(t+1)^2}}J(t),  \ \ I(t)\leq {\f{1}{(t+1)^2}}IJ(t).\label{(30)}
\end{eqnarray}
This, together with Lemma \ref{(100)}, shows
\begin{eqnarray}
\nonumber{\f {d}{dt}}J(t)&=&(2-n(\gamma-1))(t+1)E_i(t)+(t+1)^2{\f {d}{dt}}E(t)\\\nonumber
&=&(2-n(\gamma-1))(t+1)E_i(t)\\
&\leq&
\begin{cases}
   {\f{2-n(\gamma-1)}{t+1}}J(t),\qquad 1<\gamma \leq1+{\f{2}{n }},\cr
    0,\qquad\qquad\qquad \qquad \gamma >1+{\f{2}{n }},
   \end{cases}\label{(31)}
\end{eqnarray}
and
\begin{eqnarray}
\nonumber {\f {d}{dt}}I J(t)&=&(2-n(\gamma-1))(t+1)I(t)+(t+1)^2{\f {d}{dt}}IE(t)\\\nonumber
&\leq&(2-n(\gamma-1))(t+1)I(t)\\
&\leq&
 \begin{cases}
   {\f{2-n(\gamma-1)}{t+1}}I J(t),\qquad 1<\gamma \leq1+{\f{2}{n }},\cr
    0,\qquad\qquad\qquad \qquad \gamma >1+{\f{2}{n }}.
   \end{cases}\label{(32)}
\end{eqnarray}

The proof of the proposition is finished.

 It follows from  Lemma
\ref{(105)} that

\begin{Proposition}\label{(106)}
For(CNS) and (ICNS), the following estimates hold:

\begin{eqnarray}
E_i(t)\leq
\begin{cases}
   {\f{C_4}{(t+1)^{n(\gamma-1)}}},\qquad 1<\gamma \leq1+{\f{2}{n }},\cr
    {\f{C_4}{(t+1)^2}},\qquad \qquad \gamma >1+{\f{2}{n }},
   \end{cases}\label{(33)}
\end{eqnarray}
and

\begin{eqnarray}
I(t)\leq
\begin{cases}
   {\f{C_5}{(t+1)^{n(\gamma-1)}}},\qquad 1<\gamma \leq1+{\f{2}{n }},\cr
    {\f{C_5}{(t+1)^2}},\qquad \qquad \gamma >1+{\f{2}{n }},
   \end{cases}\label{(34)}
\end{eqnarray}
respectively, where $
C_4=J(0),C_5=IJ(0)$.
\end{Proposition}

Now  we are ready to prove Theorem \ref{(1000)} and Theorem
\ref{(1001)}.

 \textbf{\textbf{Proof of Theorem \ref{(1000)}}.}
Suppose that the life span of the classical solution $t=+\infty$.
Then by Proposition \ref{(104)} and Proposition \ref{(106)},  if
$1<\gamma \leq1+{\f{2}{n }}$,  we have
\begin{equation}
{\f {C_2}{G(t)^{\f {n (\gamma-1)}{2}}}}\leq
E_i(t)\leq{\f{C_4}{(t+1)^{n(\gamma-1)}}},\label{(43.01)}
\end{equation}
for all $t\geq0$.
In view of \eqref{(17)}, one has
\begin{equation}
 G(t)\leq E(0)t^2+F(0)t+G(0).\label{(43.001)}
\end{equation}
Substituting  \eqref{(43.001)} to  \eqref{(43.01)} yields
\begin{equation}
{\f {C_2}{(E(0)t^2+F(0)t+G(0))^{\f {n
(\gamma-1)}{2}}}}\leq{\f{C_4}{(t+1)^{n(\gamma-1)}}}.\label{(43.002)}
\end{equation}
Let $t$ goes to infinity,  we get
\begin{equation}
{\f {E(0)^{\f {n(\gamma-1)}{2}}J(0)}{\exp ({\f{s_1}{c_\nu}})M(0)^{\f
{(n+2)\gamma-n}{2}}}}\ge({\f{\Gamma ({\f n 2}+1)}{(\pi)^{\f n
2}}})^{\gamma-1 }
{\f{1}{2^{\f{(n+2)\gamma-n}{2}}(\gamma-1)}},\label{(43.02)}
\end{equation}
which contradicts \eqref{1000-1}. Hence  if \eqref{1000-1} holds,
then there  exists  a  time $T_1^{*}<\infty$, satisfying
\begin{equation}\label{T1}
{\f {C_2}{(E(0){T_1^{*}}^2+F(0)T_1^{*}+G(0))^{\f {n
(\gamma-1)}{2}}}}>{\f{C_3}{(T_1^{*}+1)^{n(\gamma-1)}}},
\end{equation}
such that $[0,T_1^*)$ is the life span of the classical solution.
Indeed, one can solve out $T_1^{*}$ by  \eqref{T1}. The proof of the
theorem is finished.

\textbf{\textbf{Proof of Theorem} \ref{(1001)}}. The proof is
similar to Theorem \ref{(1000)} and we omit the details here.

In \cite{Xin}, Xin investigated the blow-up of the classical
solution $(\rho, u, s)\in C^1([0,T];H^m(\mathbb{R}^n))$ with
$m>[\frac{n}{2}]+2$ and $n\ge 1$ to the Cauchy problem of (CNS)
without heat conduction if the initial density has compact support,
where $s$ is the entropy of the solution. This was generalized by
Cho and Jin to the case with heat conduction in \cite{CJ}, by
Razanova
 to the case that initial data rapidly decays at far fields in
\cite{Rozanova} and recently by
 Xin and Yan  \cite{XY} to the case that initial data
has an isolated mass group.

The following is a key lemma due to Xin \cite{Xin} which says that
if the initial density has compact support, then the compact support
will keep unchanged for all time.

\begin{Lemma}[\cite{Xin}]\label{(110)}
 For the viscous compressible  Navier-Stokes equations (CNS) with
\begin{equation}
      \mu>0, 2\mu+n\lambda>0, \mathcal{K}\ge 0,\label{(40)}
\end{equation}
\end{Lemma}
Suppose that $(\rho, u, s)\in C^1([0,T];H^m(\mathbb{R}^n))$ with
$m>[\frac{n}{2}]+2$  is a classical solution of (CNS) and the
initial density has compact support. Then the support of the density
$\rho(x,t)$ will not grow in time.

 Now we give an uniform proof of
blow-up results obtained  by Xin  \cite{Xin} for $\mathcal {K}=0$
and by Cho and Jin  \cite{CJ} for $\mathcal {K}> 0$.
\begin{Proposition}\label{(112)}
 Suppose that the assumptions of Lemma \ref{(110)} and \eqref{9-25-1} hold true. Then any smooth solution to the
 Cauchy problem of (CNS) will blow up in finite time.
\end{Proposition}
\textbf{{\em Proof.}} It follows from Lemma \ref{(110)} that the
support of the density ${\textrm supp} _x{\rho(x,t)}={\textrm supp}
_x{\rho_0(x)}\triangleq D$. Adopting the arguments in \cite{Xin}, we
have $u(x,t)=0$ in the outside of the support of density. Hence the
assumptions \eqref{9-25} are satisfied. Then
\begin{eqnarray}
\nonumber G(t)&=&{\f12}\int_{{\textrm supp} _x{\rho(x,t)}} \rho |x|^2 dx={\f12}\int_D \rho |x|^2 dx\\
&\leq&{\f12}|D|^2\int_{\mathbb{R}^n} \rho dx
 ={\f12}M(0)|D|^2.\label{(41)}
\end{eqnarray}
According to Lemma \ref{(12)}, if $1<\gamma \leq1+{\f{2}{n }}$,  we
get

\begin{equation}
{\f{n(\gamma-1)}{2}}E(0)t^2+F(0)t+G(0)\leq{\f12}M(0)|D|^2,\label{(42)}
\end{equation}
and, if $\gamma >1+{\f{2}{n }}$,  we have
\begin{equation}
  E(0)t^2+F(0)t+G(0)\leq G(t)\leq {\f12}M(0)|D|^2,\label{(43)}
\end{equation}
respectively. The inequalities (\ref{(42)}) and (\ref{(43)}) imply
that the life span is finite. The proof of the proposition is
finished.

\begin{Remark}
If the temperature $\theta(x,t)=0$ in the outside of the support of
the density, the assumption \eqref{9-25-1} is satisfied
automatically.
\end{Remark}

\section{The Proof of Theorems \ref{(1002)}-\ref{(1003+)}}

 In this subsection, we will prove Theorems \ref{(1002)}-\ref{(1003+)} which is on the blow-up of the
 classical solution to the Navier-Stokes equations with density-dependent viscosities (DICNS). In the following, the notations
$M(t),\mathbb{P}(t),G(t),F(t),IE(t),IJ(t)$ are  same as in Section
2. However, $IH(t)$ should be modified as
\begin{equation}\label{DIH}
DIH(t)=2 E_k(t)+n(\gamma-1)I(t)-\int
_{\mathbb{R}^n}[h(\rho)+ng(\rho)] (\textit{div} u) dx,
\end{equation}
where $h(\rho)=\rho^\alpha, g(\rho)=(\alpha-1)\rho^\alpha$. To prove
Theorems \ref{(1002)}-\ref{(1003+)}, the key is to obtain the upper
and lower decay rates of the potential energy $I(t)$.

The following basic relationships of  the  quantities defined in
Section 2 hold true.
\begin{Lemma}\label{(115)} For (DICNS),
we have
\begin{equation}
{\f {d}{dt}}M(t)={\f {d}{dt}}\mathbb{P}(t)=0, {\f {d}{dt}}G(t)=F(t),
{\f {d}{dt}}F(t)=DIH(t),\label{(44)}
\end{equation}
and
\begin{eqnarray}
{\f {d}{dt}}IE(t)=
\begin{cases}
\displaystyle{-\alpha\int_\mathbb{R}\rho^\alpha u_x^2dx,
\quad\,\quad\,\quad\,\quad\,\quad\,\quad\,\quad\,\quad\,\quad\,n=1;}\cr\
\displaystyle{-\int_{\mathbb{R}^n}
[h(\rho)\sum_{j=1}^n(\partial_ju_j)^2+g(\rho)(
\textrm{div}u)^2}\nonumber\\
\displaystyle{+\frac{h(\rho)}{2}\sum_{i\neq
j}(\partial_ju_i)^2+h(\rho)\sum_{i>j}(\partial_ju_i)(\partial_iu_j)]dx},
\quad\ \quad\ n\geq2.\label{(45)}
\end{cases}
\end{eqnarray}
\end{Lemma}

\vspace{3mm}

The estimates of the decay rate of $I(t)$ is divided into two steps.

\vspace{3mm}

 {\it Step 1. Lower Bounds of $I(t)$}

\vspace{3mm}

Similar to Lemma \ref{(102)}, based on Lemma \ref{(115)} and Lemma
\ref{(101)}, we can  get two-sided estimates of $F(t)$ and $G(t)$.

\begin{Lemma}\label{(117)} Let $n=1$ and suppose that Assumption 2 holds.  If $\alpha\geq1$,  then there exist two positive numbers $C_{16}$ and $C_{18}$ such that
\begin{eqnarray}
\nonumber&&({\f{\mathbb{P}(0)^2}{M(0)}}-C_{18})t+C_{16}({\f{\mathbb{P}(0)^2}{2M(0)}}-IE(0))+F(0)
\leq F(t)\\
&&\leq [\max(2,\gamma-1)IE(0)+C_{18}]t-C_{16}({\f{\mathbb{P}(0)^2}{2M(0)}}-IE(0))+F(0),\label{(48)}
\end{eqnarray}
and
\begin{eqnarray}
\nonumber&&{\f12}({\f{\mathbb{P}(0)^2}{M(0)}}-C_{18})t^2+(C_{16}({\f{\mathbb{P}(0)^2}{2M(0)}}-IE(0))+F(0))t+G(0)
\leq G(t)\\
&&\leq {\f12}[\max(2,\gamma-1)IE(0)+C_{18}]t^2+[-C_{16}({\f{\mathbb{P}(0)^2}{2M(0)}}-IE(0))+F(0)]t+G(0),\label{(49)}
\end{eqnarray}
where $C_{18}$ is any positive number between $0$ and
$\f{\mathbb{P}(0)^2}{2M(0)}$ satisfying
$C_{16}C_{18}={\f{\alpha}{4}}M(0)C_{14}^{\alpha-1}$.
\end{Lemma}
\textbf{{\em Proof.}} Using the  Cauchy's inequality, we have
\begin{eqnarray}
    \nonumber |\int _{\mathbb{R}}[h(\rho)+g(\rho)] \textit{div}u dx|
    &=&\alpha|\int_\mathbb{R}\rho^\alpha u_xdx|
     \leq C_{15}\alpha\int_\mathbb{R}\rho^\alpha dx+C_{16}\alpha\int_\mathbb{R}\rho^\alpha u_x^2 dx\\
     &\leq &C_{17} \int_\mathbb{R}\rho dx-C_{16}{\f {d}{dt}}IE(t)
     = C_{18}-C_{16}{\f {d}{dt}}IE(t),\label{(50)}
\end{eqnarray}
where $C_{16}C_{18}={\f{\alpha}{4}}M(0)C_{14}^{\alpha-1}$ and
$C_{14}=\max \rho$. Choose $C_{16}$ large enough  such that
$0<C_{18}<{\f{\mathbb{P}(0)^2}{M(0)}}$. From Lemma \ref{(115)}, we
have
\begin{eqnarray}
     \nonumber{\f {d^2}{dt^2}}G(t)
    ={\f {d}{dt}}F(t)=DIH(t)
     &\geq& 2E_k(t)+(\gamma -1)I(t)-C_{18}+C_{16}{\f {d}{dt}}IE(t)\\
    &\geq& {\f{\mathbb{P}(0)^2}{M(0)}}-C_{18}+C_{16}{\f {d}{dt}}IE(t),\label{(51)}
\end{eqnarray}
and
\begin{eqnarray}
    \nonumber {\f {d^2}{dt^2}}G(t)\leq{\f {d}{dt}}F(t)=DIH(t)
     &\le&2E_k(t)+(\gamma -1)I(t)+C_{18}-C_{16}{\f {d}{dt}}IE(t)\\
     &\leq& \max(2,\gamma-1)IE(0)+C_{18}-C_{16}{\f {d}{dt}}IE(t).\label{(52)}
\end{eqnarray}
Integrating \eqref{(51)} and \eqref{(52)} with respect to $t$ over
$[0,t]$, we finish the proof of the lemma.
\begin{Lemma}\label{(117.5)} Let $n\geq1$ and suppose that Assumption 1 hold.  If $1\leq\alpha\leq\gamma$, then there exist two positive numbers $C_{21}$ and $C_{22}$ such that
\begin{eqnarray}
\nonumber&&({\f{\mathbb{P}(0)^2}{M(0)}}-C_{22})t+C_{21}({\f{\mathbb{P}(0)^2}{2M(0)}}-IE(0))+F(0)
\leq F(t)\\
&&\leq
[\max(2,n(\gamma-1))IE(0)+C_{22}]t-C_{21}({\f{\mathbb{P}(0)^2}{2M(0)}}-IE(0))+F(0),\label{(48)}
\end{eqnarray}
and
\begin{eqnarray}
\nonumber&&{\f12}({\f{\mathbb{P}(0)^2}{M(0)}}-C_{22})t^2+(C_{21}({\f{\mathbb{P}(0)^2}{2M(0)}}-IE(0))t+G(0)
\leq G(t)\\
&&\leq
{\f12}[\max(2,n(\gamma-1))IE(0)+C_{22}]t^2+[-C_{21}({\f{\mathbb{P}(0)^2}{2M(0)}}-IE(0))+F(0)]t+G(0),\label{(49)}
\end{eqnarray}
where where $C_{22}$ is any positive number between $0$ and
$\f{\mathbb{P}(0)^2}{2M(0)}$ satisfying
$C_{21}C_{22}={\f14}n[1+n(\alpha-1)]^2(\gamma-1)^{\f{\alpha-1}{\gamma-1}}M(0)^{\f{\gamma-\alpha}{\gamma-1}}C_{13}^{\f{\alpha-1}{\gamma-1}}$.
\end{Lemma}
\textbf{{\em
Proof.}} By interpolation inequality, we have
\begin{eqnarray}
  \int_{\mathbb{R}^n} \rho^\alpha dx
  \leq (\int_{\mathbb{R}^n} \rho dx)^{\f{\gamma-\alpha}{\gamma-1}}(\int_{\mathbb{R}^n} \rho^\gamma dx)^{\f{\alpha-1}{\gamma-1}}
  \leq (\gamma-1)^{\f{\alpha-1}{\gamma-1}}M(0)^{\f{\gamma-\alpha}{\gamma-1}}C_{13}^{\f{\alpha-1}{\gamma-1}}
  \triangleq C_{19}.\label{(49.1)}
\end{eqnarray}
Here we used the Assumption 1 which means that the energy is finite.
Using Young's inequality yields
\begin{eqnarray}
|\int_{\mathbb{R}^n}[h(\rho)+n g(\rho)] (\textit{div} u) dx|
&=&[1+n(\alpha-1)]|\int_{\mathbb{R}^n} \rho^\alpha (\textit{div} u)
dx|\\\nonumber &\leq&
n^{\f{1}{2}}[1+n(\alpha-1)](\int_{\mathbb{R}^n} \rho^\alpha
dx)^{\f12}(\int_{\mathbb{R}^n} \rho^\alpha |\nabla u|^2
dx)^{\f12}\\\nonumber
&\leq& C_{20}\int_{\mathbb{R}^n} \rho^\alpha dx+C_{21}\int_{\mathbb{R}^n} \rho^\alpha |\nabla u|^2 dx\\
&\leq& C_{22}-C_{21}{\f {d}{dt}}IE(t),\label{(49.2)}
\end{eqnarray}
where $C_{21}C_{22}={\f14}C_{19}n[1+n(\alpha-1)]^2$ and we choose
$C_{21}$ is large enough such that
$C_{22}<{\f{\mathbb{P}(0)^2}{M(0)}}$. Similarly to \eqref{(51)} and
\eqref{(52)}, one can obtain
\begin{eqnarray}
   {\f {d^2}{dt^2}}G(t)={\f {d}{dt}}F(t)
    \geq {\f{\mathbb{P}(0)^2}{M(0)}}-C_{22}+C_{21}{\f {d}{dt}}IE(t),\label{(49.3)}
\end{eqnarray}
and
\begin{eqnarray}
    {\f {d^2}{dt^2}}G(t)={\f {d}{dt}}F(t)
   \leq \max(2,n(\gamma-1))E(0)+C_{22}-C_{21}{\f {d}{dt}}IE(t).\label{(49.4)}
\end{eqnarray}
Integrating \eqref{(49.3)} and \eqref{(49.4)} with respect to $t$
over $[0,t]$, we finish the proof of the lemma.

On the other hand, it follows from Lemma \ref{(103)} that

\begin{Proposition}\label{(118)} For (DICNS), we have
 \begin{equation}
  I(t)\geq {\f {C_{23}}{G(t)^{\f {n(\gamma-1)}{2}}}}.\label{(53)}
\end{equation}
where $C_{23}=C_3$.
\end{Proposition}

Hence Lemma \ref{(117)}, Lemma \ref{(117.5)} and Proposition
\ref{(118)} give the lower bounds of $I(t)$. In the next step, we
will give the upper bound of $I(t)$.

 \vspace{3mm}

 {\it Step 2. Upper Bound of $I(t)$}

\vspace{3mm}
 \begin{Lemma}\label{(119)} Let $n=1$.  Then we have
the following inequality:
\begin{eqnarray}
{\f {d}{dt}}IJ(t)+\alpha(t+1)^2\int _\mathbb{R}\rho^\alpha u_x^2 dx
\leq
\begin{cases}
{\f{3-\gamma}{t+1}}IJ(t)+\alpha(t+1)\int _\mathbb{R}\rho^\alpha u_x dx,\qquad 1<\gamma \leq3,\cr
  \alpha(t+1)\int _\mathbb{R}\rho^\alpha u_xdx,\qquad \qquad\qquad \qquad \gamma >3.
\end{cases}\label{(54)}
\end{eqnarray}
\end{Lemma}
\textbf{{\em
Proof.}} By Lemma \ref{(115)}, one can compute
\begin{eqnarray}
{\f {d}{dt}}IJ(t) &= &(3-\gamma)(t+1)I(t)+\alpha(t+1)\int
_\mathbb{R}\rho^\alpha u_x
dx-\alpha(t+1)^2\int_\mathbb{R}\rho^\alpha u_x^2dx.\label{(55)}
\end{eqnarray}
The case of $\gamma>3$  holds obviously.  The other case of
$1<\gamma \leq 3$ follows from \eqref{(30)}.

Similar to Lemma \ref{(119)}, we can get
\begin{Lemma}\label{(121)} Let $n\geq2$. Then we have
the following inequality:
\begin{eqnarray}
 {\f {d}{dt}}IJ(t) +(t+1)^2\{\int_{\mathbb{R}^n}
[\rho^\alpha\sum_{j=1}^n(\partial_ju_j)^2+(\alpha-1)\rho^\alpha(
\textrm{div}u)^2 \nonumber\\
+\frac{\rho^\alpha}{2}\sum_{i\neq
j}(\partial_ju_i)^2+\rho^\alpha\sum_{i>j}(\partial_ju_i)(\partial_iu_j)] dx\}\nonumber\\
\leq
\begin{cases}
\displaystyle{\f{2-n(\gamma-1)}{t+1}}IJ(t)+[1+n(\alpha-1)](t+1)\int_{\mathbb{R}^n}
\rho^\alpha (\textrm{div}u) dx,\qquad 1<\gamma \leq1+{\f{2}{n }},\cr
\displaystyle{[1+n(\alpha-1)](t+1)\int_{\mathbb{R}^n} \rho^\alpha
(\textrm {div}u) dx,}\qquad \qquad \qquad\qquad \qquad \gamma
>1+{\f{2}{n }}.
\end{cases}\label{(61)}
\end{eqnarray}
\end{Lemma}

In view of Lemma \ref{(119)} and Lemma \ref{(121)}, we can  obtain
the upper bound of $I(t)$.
\begin{Proposition}\label{(120)} Let $n=1$.  If $\alpha\geq\gamma$,  we have
the following estimates:
\begin{eqnarray}
I(t)\leq
\begin{cases}
    C_{26}(1+t)^{1-\gamma}\exp(-{\f{C_{25}}{t+1}}),\qquad\qquad 1<\gamma \leq3,\cr
    C_{26}(1+t)^{-2}\exp(-{\f{C_{25}}{t+1}}), \qquad  \qquad \qquad \gamma >3.\label{(57)}
\end{cases}
\end{eqnarray}
where $C_{25}={\f{\alpha(\gamma-1)}{4}}C_{14}^{\alpha-\gamma}$ and $ C_{26}=IJ(0)\exp({\f{\alpha(\gamma-1)}{4}}C_{14}^{\alpha-\gamma})$.
\end{Proposition}
\textbf{{\em Proof.}} It follows from Young's  inequality and
\eqref{(30)} that
\begin{eqnarray}
\nonumber \alpha(t+1)\int _\mathbb{R}\rho^\alpha u_x dx
&\leq& {\f{\alpha}{4}}\int _\mathbb{R}\rho^\alpha dx+\alpha(t+1)^2\int _\mathbb{R}\rho^\alpha u_x^2 dx\\\nonumber
&\leq& C_{24}\int _\mathbb{R}\rho^\gamma dx+\alpha(t+1)^2\int _\mathbb{R}\rho^\alpha u_x^2 dx\\
&\leq& {\f{C_{25}}{(t+1)^2}}IJ(t)+\alpha(t+1)^2\int
_\mathbb{R}\rho^\alpha u_x^2 dx.\label{(58)}
\end{eqnarray}
Here we used the fact that $\alpha\geq\gamma$. Substituting
\eqref{(58)} into \eqref{(54)} yields
\begin{eqnarray}
{\f {d}{dt}}IJ(t) \leq
\begin{cases}
    {\f{C_{25}}{(t+1)^2}}IJ(t)+{\f{3-\gamma}{t+1}}IJ(t),\qquad 1<\gamma \leq3,\cr
    {\f{C_{25}}{(t+1)^2}}IJ(t),\qquad \qquad \qquad \qquad \gamma >3. \label{(59)}
\end{cases}
\end{eqnarray}
Using Gronwall's inequality, we obtain
 \begin{eqnarray}
IJ(t)\leq
\begin{cases}
    C_{26}(1+t)^{3-\gamma}\exp(-{\f{C_{25}}{t+1}}),\qquad\qquad 1<\gamma \leq3,\cr
    C_{26}\exp(-{\f{C_{25}}{t+1}}), \qquad \qquad \qquad  \qquad \qquad \gamma >3.
\end{cases}\label{(60)}
\end{eqnarray}
This, together with \eqref{(30)}, shows
\begin{eqnarray}
I(t)\leq
\begin{cases}
    C_{26}(1+t)^{1-\gamma}\exp(-{\f{C_{25}}{t+1}}),\qquad\qquad 1<\gamma \leq3,\cr
    C_{26}(1+t)^{-2}\exp(-{\f{C_{25}}{t+1}}), \qquad  \qquad \qquad \gamma >3.
\end{cases}
\end{eqnarray}
Here  $C_{24}={\f{\alpha}{4}}C_{14}^{\alpha-\gamma},
C_{25}={\f{\alpha(\gamma-1)}{4}}C_{14}^{\alpha-\gamma}$ and $
C_{26}=IJ(0)\exp({\f{\alpha(\gamma-1)}{4}}C_{14}^{\alpha-\gamma})$.
The proof of the proposition is finished.

\begin{Proposition}\label{(120.5)} Let $n\ge 1$.  If $\frac{\gamma+1}{2}<\alpha\le \gamma$ and $1<\gamma<1+\frac2n$,  then
we have
the following estimates:
\begin{eqnarray}
I(t)\leq
\begin{cases}
C_{29}(1+t)^{-n(\gamma-1)}\exp(-{\f{C_{27}}{(t+1)^{C_{28}}}})+\frac{1}{(1+t)^2},
\quad \quad \quad \quad\quad \quad    \ n\geq2, \cr\
C_{31}(1+t)^{1-\gamma}\exp(-{\f{C_{30}}{(t+1)^{C_{28}}}})+\frac{1}{(1+t)^2},
\quad \quad \quad \quad \quad \quad \quad    \ n=1.
\end{cases}\label{(66)}
\end{eqnarray}
Here
\begin{align}
\nonumber&C_{27}=\exp({\f{\gamma-1}{4(\alpha-1)(2\alpha-\gamma-1)}}[1+n(\alpha-1)]^2(\gamma-1)^{\f{\alpha-1}{\gamma-1}})
M(0)^{\f{\gamma-\alpha}{\gamma-1}});\\\nonumber
&C_{28}={\f{2\alpha-\gamma-1}{\gamma-1}};\\\nonumber &C_{29}=IJ(0)
\exp({\f{\gamma-1}{4(\alpha-1)(2\alpha-\gamma-1)}}[1+n(\alpha-1)]^2
(\gamma-1)^{\f{\alpha-1}{\gamma-1}}M(0)^{\f{\gamma-\alpha}{\gamma-1}});\\\nonumber
\nonumber
&C_{30}=\exp({\f{\alpha(\gamma-1)}{4(2\alpha-\gamma-1)}}(\gamma-1)^{\f{\alpha-1}{\gamma-1}}
M(0)^{\f{\gamma-\alpha}{\gamma-1}});\\\nonumber &C_{32}=IJ(0)
\exp({\f{\alpha(\gamma-1)}{4(2\alpha-\gamma-1)}}(\gamma-1)^{\f{\alpha-1}{\gamma-1}}
M(0)^{\f{\gamma-\alpha}{\gamma-1}}).
\end{align}
\end{Proposition}

\textbf{{\em Proof of Proposition \ref{(120.5)}.}}  It follows from
\eqref{(49.1)} that
\begin{eqnarray}
  \int_{\mathbb{R}^n} \rho^\alpha dx
  \leq (\int_{\mathbb{R}^n} \rho dx)^{\f{\gamma-\alpha}{\gamma-1}}(\int_{\mathbb{R}^n} \rho^\gamma dx)^{\f{\alpha-1}{\gamma-1}}
  \leq (\gamma-1)^{\f{\alpha-1}{\gamma-1}}M(0)^{\f{\gamma-\alpha}{\gamma-1}}(t+1)^{-{\f{2(\alpha-1)}{\gamma-1}}}
  IJ(t)^{\f{\alpha-1}{\gamma-1}}
 \label{(62)}
\end{eqnarray}
for $1\le \alpha\le \gamma$. Consequently, when  $n\geq2$, one has
\begin{eqnarray}
 \nonumber &&[1+n(\alpha-1)](t+1)\int_{\mathbb{R}^n} \rho^\alpha \textrm{div}u dx \\ \nonumber
  &\leq&{\f{1}{4(\alpha-1)}}[1+n(\alpha-1)]^2\int_{\mathbb{R}^n} \rho^\alpha dx+(\alpha-1)(t+1)^2\int_{\mathbb{R}^n} \rho^\alpha\textrm({div}u)^2dx \\\nonumber
  &\leq&{\f{1}{4(\alpha-1)}}[1+n(\alpha-1)]^2(\gamma-1)^{\f{\alpha-1}{\gamma-1}}
  M(0)^{\f{\gamma-\alpha}{\gamma-1}}(t+1)^{-{\f{2(\alpha-1)}{\gamma-1}}}
  IJ(t)^{\f{\alpha-1}{\gamma-1}}\\
  &&+(\alpha-1)(t+1)^2\int_{\mathbb{R}^n} \rho^\alpha\textrm({div}u)^2dx. \ \label{(62.1)}
\end{eqnarray}
When $n=1$, one has
\begin{eqnarray}
\nonumber &&\alpha(t+1)\int _\mathbb{R}\rho^\alpha u_x dx \leq
{\f{\alpha}{4}}\int _\mathbb{R}\rho^\alpha dx+\alpha(t+1)^2\int
_\mathbb{R}\rho^\alpha u_x^2 dx\\\nonumber &\leq&
{\f{\alpha}{4}}(\gamma-1)^{\f{\alpha-1}{\gamma-1}}
M(0)^{\f{\gamma-\alpha}{\gamma-1}}(t+1)^{-{\f{2(\alpha-1)}{\gamma-1}}}
  IJ(t)^{\f{\alpha-1}{\gamma-1}}\\
  &&+\alpha(t+1)^2\int _\mathbb{R}\rho^\alpha u_x^2 dx. \ \label{(62.2)}
\end{eqnarray}
Putting \eqref{(62.1)} and \eqref{(62.2)} into \eqref{(61)},  we
have
\begin{eqnarray}
 {\f {d}{dt}}IJ(t)
\leq
\begin{cases}
{\f{1}{4(\alpha-1)}}[1+n(\alpha-1)]^2(\gamma-1)^{\f{\alpha-1}{\gamma-1}}
 M(0)^{\f{\gamma-\alpha}{\gamma-1}}
(t+1)^{-{\f{2(\alpha-1)}{\gamma-1}}}
  IJ(t)^{\f{\alpha-1}{\gamma-1}}\\
 +{\f{2-n(\gamma-1)}{t+1}}IJ(t), \quad \quad \quad \quad \quad \quad \quad \quad \quad \quad \quad \quad \quad \quad  \quad \quad \quad \quad \ n\geq2, \cr\
 {\f{\alpha}{4}}(\gamma-1)^{\f{\alpha-1}{\gamma-1}}
 M(0)^{\f{\gamma-\alpha}{\gamma-1}}(t+1)^{-{\f{2(\alpha-1)}{\gamma-1}}}
  IJ(t)^{\f{\alpha-1}{\gamma-1}}+{\f{3-\gamma}{t+1}}IJ(t), \ n=1
 \end{cases}\label{(63)}
\end{eqnarray}
for $1<\gamma \le 1+{\f{2}{n }}$.

If $IJ(t)\leq1$, one deduces that $I(t)\leq
\frac{1}{(1+t)^2}IJ(t)\le {\f{1}{(t+1)^2}}$.

If $IJ(t)>1$, since $\alpha\le\gamma$, one has
\begin{eqnarray}
{\f {d}{dt}}IJ(t) \leq
\begin{cases}
\{{\f{2-n(\gamma-1)}{t+1}}+{\f{1}{4(\alpha-1)}}[1+n(\alpha-1)]^2(\gamma-1)^{\f{\alpha-1}{\gamma-1}}\\
\times M(0)^{\f{\gamma-\alpha}{\gamma-1}}
(t+1)^{-{\f{2(\alpha-1)}{\gamma-1}}}\}
  IJ(t), \quad \quad \quad \quad \quad \quad \quad \quad \quad  \ n\geq2, \\
 [{\f{3-\gamma}{t+1}}+{\f{\alpha}{4}}(\gamma-1)^{\f{\alpha-1}{\gamma-1}}
 M(0)^{\f{\gamma-\alpha}{\gamma-1}}(t+1)^{-{\f{2(\alpha-1)}{\gamma-1}}}] IJ(t).  \quad \quad  n=1,
 \end{cases}\label{(64)}
\end{eqnarray}
Thanks to $\gamma<1+\frac2n$  and $\frac{\gamma+1}{2}<\alpha$, we
obtain by Gronwall's inequality that
\begin{eqnarray}
IJ(t)\leq
\begin{cases}
C_{29}(1+t)^{2-n(\gamma-1)}\exp(-{\f{C_{27}}{(t+1)^{C_{28}}}}),\quad
\quad \quad
n\geq2, \\
C_{31}(1+t)^{3-\gamma}\exp(-{\f{C_{30}}{(t+1)^{C_{28}}}}), \quad
\quad \quad \quad \quad  n=1.
\end{cases}\label{(65)}
\end{eqnarray}
which implies
\begin{eqnarray}
I(t)\leq \frac{1}{(1+t)^2}IJ(t)\le
\begin{cases}
C_{29}(1+t)^{-n(\gamma-1)}\exp(-{\f{C_{27}}{(t+1)^{C_{28}}}}),
\quad \quad \quad \quad\quad \quad    \ n\geq2, \cr\
C_{31}(1+t)^{1-\gamma}\exp(-{\f{C_{30}}{(t+1)^{C_{28}}}}), \quad
\quad \quad \quad \quad \quad \quad    \ n=1.
\end{cases}\label{(66)}
\end{eqnarray}
Here
\begin{align}
\nonumber&C_{27}=\exp({\f{\gamma-1}{4(\alpha-1)(2\alpha-\gamma-1)}}[1+n(\alpha-1)]^2(\gamma-1)^{\f{\alpha-1}{\gamma-1}})
M(0)^{\f{\gamma-\alpha}{\gamma-1}});\\\nonumber
&C_{28}={\f{2\alpha-\gamma-1}{\gamma-1}};\\\nonumber &C_{29}=IJ(0)
\exp({\f{\gamma-1}{4(\alpha-1)(2\alpha-\gamma-1)}}[1+n(\alpha-1)]^2
(\gamma-1)^{\f{\alpha-1}{\gamma-1}}M(0)^{\f{\gamma-\alpha}{\gamma-1}});\\\nonumber
\nonumber
&C_{30}=\exp({\f{\alpha(\gamma-1)}{4(2\alpha-\gamma-1)}}(\gamma-1)^{\f{\alpha-1}{\gamma-1}}
M(0)^{\f{\gamma-\alpha}{\gamma-1}});\\\nonumber &C_{31}=IJ(0)
\exp({\f{\alpha(\gamma-1)}{4(2\alpha-\gamma-1)}}(\gamma-1)^{\f{\alpha-1}{\gamma-1}}
M(0)^{\f{\gamma-\alpha}{\gamma-1}}).
\end{align}

Combining the cases $IJ(t)\le 1$ and  $IJ(t)> 1$, we finish the
proof of the proposition.

\vspace{3mm}

Now, we are ready to prove Theorem \ref{(1002)} and Theorem
\ref{(1003)}. The key idea is to compare the coefficients in the
lower bounds and the upper bounds of the potential energy $I(t)$,
which is similar to the proof of Theorem \ref{(1000)} and we omit
the details here.

\textbf{{\em Proof of Theorem \ref{(1002)}.}} Theorem \ref{(1002)}
follows from Lemma \ref{(117)}, Proposition \ref{(118)} and
Proposition \ref{(120)}.

\textbf{{\em Proof of Theorem \ref{(1003)}.}}Theorem \ref{(1003)}
follows from Lemma \ref{(117.5)}, Proposition \ref{(118)} and
Proposition \ref{(120.5)}.

\textbf{{\em Proof of Theorem \ref{(1003+)}.}} Theorem \ref{(1003+)}
follows from Lemma \ref{(117.5)}, Proposition \ref{(118)} and
Proposition \ref{(120.5)}.

\end{document}